# A linear programming based approach for determining maximal closest reference set in DEA

I. Roshdi[a*], M. Mehdiloozad[b], D. Margaritis[c]

[a] *Department of Mathematics, Islamic Azad University, Semnan Branch, Semnan, Iran*

[b] *Department of Mathematics, College of Sciences, Shiraz University, Shiraz, Iran*

[c] *University of Auckland Business School, Auckland, New Zealand*

**Abstract**

Identification of the reference set for each decision making unit (DMU) is a main concern in the data envelopment analysis (DEA). All of the methods developed to date have been focused on finding the furthest reference DMUs. In this paper, we introduce the new notion of maximal closest reference set (MCRS) containing the maximum number of closest reference DMUs to the assessed DMU. Then, we develop a linear programming (LP) model for determining the MCRS for each inefficient DMU. We illustrate our method through a numerical example.

**1. Introduction**

Data envelopment analysis (DEA), initiated and developed by Charnes et al. [1], is a non-parametric LP-based technique for evaluating the relative efficiency of a set of homogeneous decision making units (DMUs). For an inefficient DMU, as well as determining its efficiency score, DEA provides a *reference set* composed of efficient DMUs to which the assessed DMU is directly compared for estimating its efficiency score. Those DMUs forming the reference set of an inefficient DMU are known as the reference DMUs or peer group for it. Form the target setting point of view, the reference DMUs can be regarded as the realistic targets of inefficient DMUs for efficiency improvement. Obviously, the more the peers are "similar" to the assessed DMU, the more the comparison is meaningful, and the less levels of operation for the inputs and outputs of inefficient DMUs are needed to make them efficient. In this sense, the reference DMUs that the smallest modifications in the inputs and outputs of the given DMU is required to reach them, are as much similar as possible to the assessed DMU.

Finding the reference set of DMUs has recently received attention in the DEA literature because of its various applications such as (i) returns to scale identification [2] (ii) sensitivity and stability analysis [3] (iii) ranking of DMUs [4] (iv) benchmarking and target setting [5, 6]. However, the methods developed to date [2, 7, 8] suffer from a common drawback: finding the furthest reference DMUs rather than the nearest ones. This drawback is a natural consequence of developing these methods based on the conventional DEA models (e.g., see the CCR and BCC models discussed in [9], which yield the furthest efficient projection points for inefficient DMUs.

The purpose of the present study is to remedy the above-mentioned drawback of the existing methods, by developing a new method for determining the maximal closest reference set for DMUs. To this purpose, we exploit the work of Aparicio et al. [10] who developed an innovative single-stage approach for finding the closest Pareto-efficient projection points for a given inefficient DMU. Based on their work, our method is presented as a two-staged approach. In the first stage, using the approach of Aparicio et al. [10], we identify a closest Pareto-efficient projection point for the assessed inefficient DMU. Having identified a closest Pareto-efficient projection point, we then solve a LP problem to determine all the closest reference DMUs associated with this projection point, as the second stage. The most significant advantage of our method is the identification of the types and amounts of inputs and outputs needed for an inefficient DMU to perform as an efficient target with the less possible effort. This advantage is very important from a managerial perspective, because it provides the decision maker with more flexible and realistic benchmarking and target setting. In addition, since our approach is based on solving LP problems, it is extremely appropriate for practical applications.

---

[*] Corresponding Author, Email: i.roshdi@gmail.com

The reminder of this paper is organized as follows. In the next section, we develop our proposed method to determine the MCRS. Section 3 provides a numerical example to illustrate the proposed approach. Finally, conclusions are stated in Section 4.

**2. Characterizing and Finding the Maximal Closest Reference Set**

In this section, we devise a new LP-based approach for determining the MCRS of each DMU through the incorporation of the SCSC (Strong Complementary Slackness Condition [9]) into the similarity concept [10]. Our study in this work is under the constant returns to scale assumption. However, the results can be recast for the variable returns to scale [11] with some simple changes.

We begin by modifying the mADD model, proposed by Aparicio et al. [10], as model (1) with replacing *SE* instead of *E*. The sets *SE* and *E* respectively denote the sets of all *efficient* and *extreme-efficient* [12] DMUs.

$$
\begin{aligned}
Min \quad & \sum_{i=1}^{m} s_i^- + \sum_{r=1}^{s} s_r^+ \\
s.t. \quad & \sum_{j \in SE} \lambda_j x_j = x_{io} - s_i^-, && i = 1,...,m, \\
& \sum_{j \in SE} \lambda_j y_j = y_{ro} + s_r^+, && r = 1,...,s, \\
& \sum_{i=1}^{m} v_i x_{ij} - \sum_{r=1}^{s} u_r y_{rj} - d_j = 0, && j \in SE, \\
& d_j \leq M b_j, && j \in SE, \\
& \lambda_j \leq M(1 - b_j), && j \in SE, \\
& \lambda_j, d_j \geq 0, \ b_j \in \{0,1\}, && j \in SE, \\
& v_i \geq 1, \ u_r \geq 1, s_i^- \geq 0, \ s_r^+ \geq 0, \ \forall i, \ \forall r.
\end{aligned} \quad (1)
$$

Let $P = (x_o^P, y_o^P) = (x_o - S^{-*}, y_o + S^{+*})$ be the closest projection point of the assessed unit, $DMU_o$, onto the efficient frontier obtained via solving model (1). Coordinates of this point give the levels of the closest Pareto-efficient target dominating $DMU_o$ [10]. However, the obtained projection point may be an unobserved (virtual) DMU and cannot be regarded as a benchmark for $DMU_o$. In such situation, one could identify the reference DMUs that the projection point P can be expressed as a positive combination of them, and then compare the performance of $DMU_o$ with the reference DMUs as its benchmarks. This comparison will be more valid and reliable when the identified reference DMUs are closer to the assessed DMU. Therefore, we need to characterize the set containing all the closest reference DMUs. To do this, first, we present some basic definitions based on the optimal solutions of model (1).

**Definition 1.** (*Closest-Reference Set*) Let $(\lambda^*, S^{-*}, S^{+*}, d^*, U^*, V^*)$ be an optimal solution of (1) corresponding to a given DMU$_o$ and $P = (x_o^P, y_o^P)$ be the associated projection point. We define the set of all DMUs corresponding to positive $\lambda_j^*$ as the closest reference set (CRF) or the closest peer group to DMU$_o$ and denote it by $R_{op}^C$ i.e.,

$$R_{op}^C = \{DMU_j \mid \lambda_j^* > 0 \text{ in some optimal solution of (1)}\} \subseteq \{1,...,n\}. \tag{2}$$

Each member of the set $R_{op}^C$ is also called a *closest reference DMU* to DMU$_o$.

Because of the occurrence of multiple optimal solutions in model (1), multiple projection points may be obtained for it. Therefore, after this, we suppose that the point $P = (x_o^P, y_o^P)$ is a fixed projection of DMU$_o$ obtained via solving (1). Corresponding to this projection point, multiple optimal solutions may occur for the vector $\lambda$ and the CRFs will not be unique, accordingly. To formulate this, we introduce the following definition.

**Definition 2.** (*Maximal Closest Reference Set*) We define the set containing *all* the CRFs for DMU$_o$ associated to the projection point P as the maximal closest reference set (MCRF) for DMU$_o$, and denote it by $R_{op}^{MC}$.

From the above definition, the set $R_{op}^{MC}$ contains the maximum number of the closest reference DMUs of DMU$_o$ associated with the projection point P. Thus, our goal is concentrated on the characterization of this set. To do this, we define $E_o^P$ as the set of indices of all the efficient DMUs on the supporting hyperplane passing through the projection point P or formally as:

$$E_o^P = \{j \in SE \mid d_j^* = 0 \text{ in (1)}\}. \tag{3}$$

Then, by means of the following lemma, we prove that all members of $R_{op}^{MC}$ are efficient and lie on the hyperplane binding at point P.

**Lemma 1.** Let $(\lambda^*, S^{-*}, S^{+*}, d^*, U^*, V^*)$ be an optimal solution of (1) in evaluating DMU$_o$. Then, all members of $R_{op}^{MC}$ lie on the hyperplane $U^*y - V^*x = 0$, i.e., $R_{op}^{MC} \subseteq \{DMU_j \mid j \in E_o^P\}$.

***Proof.*** Since the hyperplane $U^*y - V^*x = 0$ passes through the point $P = (x_o^P, y_o^P)$, it also passes through each DMU$_j$ that enters actively- i.e., with positive $\lambda_j^* > 0$ in a combination as

$$\sum_{j \in SE} \lambda_j^* x_j = x_o^P, \quad \sum_{j \in SE} \lambda_j^* y_j = y_o^P.$$

Thus, for each DMU$_j$ belonging to $R_{op}^{MC}$, we conclude that $U^*y_j - V^*x_j = 0$. This completes the proof. ∎

Now, we formulate the following LP problem to find the MCRF for inefficient DMU$_o$:

$$\begin{aligned}
Max \quad & \eta \\
s.t. \quad & \sum_{j \in E_o^P} \mu_j x_{ij} = x_{io}^P, & i = 1,...,m, \\
& \sum_{j \in E_o^P} \mu_j y_{rj} = y_{ro}^P, & r = 1,...,s, \\
& \sum_{r=1}^{s} u_r y_{ro}^P - \sum_{i=1}^{m} v_i x_{io}^P = 0, \\
& \sum_{r=1}^{s} u_r y_{rj} - \sum_{i=1}^{m} v_i x_{ij} + t_j = 0, \quad j \in E_o^P, \\
& v_i \geq 1, \ u_r \geq 1, & \forall i, \forall r, \\
& \mu_j + t_j \geq \eta & j \in E_o^P, \\
& \mu_j, t_j \geq 0, & j \in E_o^P.
\end{aligned} \quad (4)$$

As we will demonstrate via the following lemma, the constraints of this model are the set of constraints of a pair of primal and dual LP problems, as well as two other constraints imposing the optimality condition and the SCSC on the variables $\mu_j$ and $d_j$ as a pair of primal and dual complementary variables.

**Lemma 2.** The optimal objective of model (4) is positive i.e., $\eta^* > 0$.

*Proof.* Consider the following pair of primal and dual models:

| **Primal** | **Dual** |
|---|---|
| $Min \ -\sum_{i=1}^{m} s_i^- - \sum_{r=1}^{s} s_r^+$ | $Max \ \sum_{r=1}^{s} u_r y_{ro}^P - \sum_{i=1}^{m} v_i x_{io}^P$ |
| $s.t. \ \sum_{j \in E_o^P} \mu_j x_{ij} = x_{io}^P - s_i^-, \ i=1,...,m,$ | $s.t. \ \sum_{r=1}^{s} u_r y_{rj} - \sum_{i=1}^{m} v_i x_{ij} + t_j = 0, \ j \in E_o^P,$ |
| $\sum_{j \in E_o^P} \mu_j y_{rj} = y_{ro}^P + s_r^+, \ r=1,...,s,$ | $v_i \geq 1, \ u_r \geq 1, \quad \forall i, \forall r.$ |
| $\mu_j \geq 0, \ s_i^- \geq 0, \ s_r^+ \geq 0, \ \forall j, \forall i, \forall r.$ | |
| (5) | (6) |

Because the projection P is a Pareto-efficient point, the primal model (5) has a finite optimal objective equal to zero. From the duality theorems in linear programming, it then follows that the dual model (6) is also feasible with the optimal objective value of zero. Thus, from the SCSC theorem, there exists a pair of optimal solutions of (5) and (6), $(\mu^*, S^{-*} = 0, S^{+*} = 0)$ and $(t^*, U^*, V^*)$, such that the SCSC is held between them. Now, since $(\mu^*, t^*, U^*, V^*, \eta^*)$ is a feasible solution of (4) with $\eta = Min_{j \in E_o^P} \{\mu_j^* + t_j^*\} > 0$ and model (4) is to be maximized, so $\eta^* > 0$. ■

As a corollary of the above-proved lemmas, the following theorem holds.

**Theorem 1.** Suppose that $(\mu^*, t^*, U^*, V^*, \eta^*)$ is an optimal solution of model (4). Then, $R_{op}^{MC} \subseteq \{DMU_j | \mu_j^* > 0\}$.

To sum up the characterization completed in the above examination, the proposed approach for determining the MCRS of DMU$_o$ can be summarized as the following two-step manner:

*Step 1.* Evaluate DMU$_o$ via the MILP problem (1) and identify the set $E_o^P$.

*Step 2.* Evaluate the projection point obtained in Step 1 via the LP problem (4).

## 3. A Numerical Example

We now provide a numerical example to illustrate the application of our proposed approach. The sample consists of nine DMUs, each using two inputs to produce the same amount of output. The input-output data are recorded in Table 1.

**Table 1** Data set for the numerical example

|       | DMU1 | DMU2 | DMU3 | DMU4 | DMU5 | DMU6 | DMU7 | DMU8 | DMU9 |
|-------|------|------|------|------|------|------|------|------|------|
| $x_1$ | 1    | 2    | 4    | 6    | 3    | 3    | 7    | 5    | 9    |
| $x_2$ | 7    | 5    | 3    | 2    | 4    | 8    | 4    | 3    | 3    |
| $y$   | 1    | 1    | 1    | 1    | 1    | 1    | 1    | 1    | 1    |

Fig. 1 illustrates the data listed in Table 1, in which the DMUs 1, 2, 3, 4 and 5 are actually extreme efficient (expect DMU 5) and the remaining four units are inefficient.

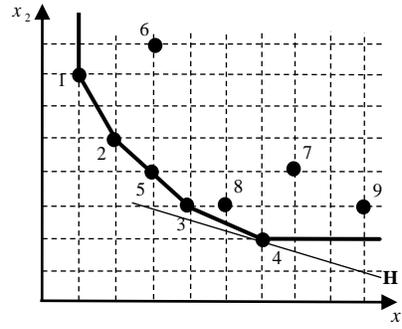

Fig. 1 Graphical representation of the data given in Table 1

To examine the application of the approach discussed in Section 2, we, specifically, have reported the furthest and closest projection points for all inefficient DMUs respectively obtained from the additive model [13] and model (1) in Table 2. Apart from the projection points given in Table 2, the maximal furthest reference set (FRS) and closest reference set (CRS) associated with them are also available.

**Table 2** Furthest and closest projection points selected for each inefficient DMU and their associated maximal FRSs and CRFs

|      | Furthest         |                    | Closest              |                 |
|------|------------------|--------------------|----------------------|-----------------|
|      | Projection Point | Maximal FRS        | Projection Point     | Maximal CRS     |
| **DMU6** | DMU5         | {DMU2,DMU3,DMU5}   | (3.0000,8.0000,1.5556) | {DMU1,DMU2}   |
| **DMU7** | DMU5         | {DMU2,DMU3,DMU5}   | (7.0000,2.3333,1.1667) | {DMU4}        |
| **DMU8** | DMU3         | {DMU3}             | (5.0000,3.0000,1.1000) | {DMU3,DMU4}   |
| **DMU9** | DMU3         | {DMU3}             | (9.0000,3.0000,1.5000) | {DMU4}        |

As can be seen from Table 2, there exist significant differences between the furthest and closest targets computed via the additive model and model (1) for the inefficient DMUs. Consequently, the obtained maximal FRSs are different from the maximal CRSs. For example, the only closest reference DMU of DMU7 is DMU4, while the furthest references for it are the DMUs 2, 3 and 5. These differences are of crucial importance in benchmarking and target setting. For instance, the inefficient DMU6 should decrease at least 4 and 1 unit in the first and second inputs to reach each member of its FRS, {DMU2, DMU3, DMU5}; while, it needs 1 and 2 unit improvements in its first and second inputs to reach the closest reference DMU1 belonging its CRS, {DMU1, DMU2}.

To highlight the role of model (4) in determining the MCRS, consider, for instance, the case of DMU9: In the evaluation of DMU9, model (1) has alternative optimal solutions, each corresponds to a supporting hyperplanes passing through the projection point $(9,3,1.5) = 1.5 \times (6,2,1) = 1.5 \times \text{DMU4}$, such as H and the hyperplane passing through the line segment joining DMU3 and DMU4 (see Fig. 1). From the results of Table 2, DMU3, which lies on the latter hyperplane, does not belong to the MCRS of DMU9 (i.e., $\lambda^*_{\text{DMU3}} = 0$ and $t^*_{\text{DMU3}} > 0$ in the optimal solution of model (4)). On the other hand, DMU4 is the only reference DMU for DMU9 that is on the hyperplane H. Therefore, from the optimal solutions to (1), model (4) select those associated with the hyperplanes that exactly all the members of MCRS lie on them.

**4. Conclusions**

In this paper, we introduced the useful notion of the MCRS containing the maximum number of the closest reference DMUs to the assessed DMU. Then, we proposed an LP-based method for determining the MCRS. We based our method on the approach proposed in [10] for identifying the closest Pareto-efficient targets of inefficient DMUs. Having determined an efficient target to the assessed DMU, the proposed procedure involves solving a unique LP problem.

Finding the references of DMUs is useful in sensitivity and stability analysis, the status of returns to scale of a DMU, ranking of DMUs. Specifically, proposing a way with the less possible adjustments of inputs and outputs, it is very useful in benchmarking and target setting. In fact, the coordinates of a closest reference DMU can be interpreted as the observed "target" levels of operation of inputs and outputs that indicate how the evaluated DMU can be improved to perform efficiently. Thus, the more the references are closer to the under assessment unit, the more they are preferable and meaningful. Based on this finding, the extension of the proposed method in the application mentioned above is suggested as future research directions.